\RequirePackage{fix-cm}
\documentclass[smallextended]{svjour3}       
\smartqed  
\usepackage[parfill]{parskip}    
\usepackage{amsmath, amssymb, latexsym, enumerate, mathrsfs}
\usepackage{graphicx}
\usepackage{epstopdf}
\usepackage{subfig}

\usepackage{hyperref}
\usepackage{array}
\usepackage{times}
\usepackage{moreverb}
\usepackage{marvosym}
\DeclareMathOperator*{\arginf}{arg\,inf}

\DeclareMathOperator*{\essinf}{ess\,inf}

\DeclareMathOperator{\blockdiag}{block\,diag}
\newtheorem{assumption}[lemma]{Assumption}
\newcommand{\etal}{\mbox{et al.}}

\begin{document}

\title{Acceptable risks and related decision problems with multiple risk-averse agents
}
\subtitle{}


\author{Getachew K. Befekadu \and Eduardo L. Pasiliao
}


\institute{G. K. Befekadu \at
	          National Research Council, Air Force Research Laboratory \& Department of Industrial System Engineering, University of Florida - REEF, 1350 N. Poquito Rd, Shalimar, FL 32579, USA.
	          \and
	          E. L. Pasiliao \at
	          Munitions Directorate, Air Force Research Laboratory, 101 West Eglin Blvd, Eglin AFB, FL 32542, USA. \\
}

\date{Received: November 15, 2016 / Accepted: date}

\maketitle

\begin{abstract}
In this paper, we consider a risk-averse decision problem for controlled-diffusion processes, with dynamic risk measures, in which multiple risk-averse agents choose their decisions in such a way to minimize their individual accumulated risk-costs over a finite-time horizon. In particular, we introduce multi-structure dynamic risk measures induced from conditional $g$-expectations, where the latter are associated with the generator functionals of certain BSDEs that implicitly take into account the risk-cost functionals of the risk-averse agents. Here, we also require that such solutions of the BSDEs to satisfy a stochastic viability property with respect to a given closed convex set. Moreover, using a result similar to that of the Arrow-Barankin-Blackwell theorem, we establish the existence of consistent optimal decisions for the risk-averse agents, when the set of all Pareto optimal solutions, in the sense of viscosity, for the associated dynamic programming equations is dense in the given closed convex set. Finally, we briefly comment on the characteristics of acceptable risks vis-\'{a}-vis some uncertain future costs or outcomes, in which results from the dynamic risk analysis constitute part of the information used in the risk-averse decision criteria.
\end{abstract}

\section{Introduction}	 \label{S1}
During the past decades, since the early work of Borch \cite{Bor62}, there have been studies clarifying appropriate solution concepts, such as equilibrium solutions, in connection with optimal risk allocations and risk aversions that are almost exclusively in the context of mathematical economics and insurance (e.g., see \cite{BarElK05}, \cite{BuhJ79}, \cite{ChatDT00}, \cite{DanS07}, \cite{Den01}, \cite{Gerb78} and \cite[pp.~88-96]{Gerb79} and the references therein). On the other hand, since the seminal work of Artzner \etal \,\cite{ArtDEH99}, some interesting studies on characterizing axiomatic concepts such as dynamic risk measures, coherency, consistency and convexity have been reported in the literature (e.g., see \cite{DetS05}, \cite{Pen04}, \cite{Ros06}, \cite{FolS02} or \cite{CorHMP02} and the references therein). More importantly, we observe that the contributions of these studies at least focused on three interrelated questions: (i) the first one is a conceptual or purely modeling question that deals with, for example, how should solution concepts -- such as equilibrium solutions in the context of optimal risk allocations and risk aversions -- be defined given new decision theoretic foundations; (ii) the second one is a general question on the behavioral implication -- associated with the consistency and coherency of ``rational" decisions -- of such newly introduced solution concepts; and (iii) the last one is an insight question -- where such solution concepts or the associated innovative results might bring to applied contexts in mathematical economics, finance, engineering and elsewhere.

In this paper, we consider a risk-averse decision problem for controlled-diffusion processes, with dynamic risk measures, in which multiple risk-averse agents choose their decisions in such a way to minimize their individual accumulated risk-costs over a finite-time horizon. We specifically introduce multi-structure dynamic risk measures induced from conditional $g$-expectations, where the latter are associated with the generator functionals of certain BSDEs that implicitly take into account the risk-cost functionals of the risk-averse agents. Here, we also require that such solutions of the BSDEs to satisfy a stochastic viability property with respect to a given closed convex set. Further, using a result similar to the Arrow-Barankin-Blackwell theorem, we establish the existence of consistent optimal decisions for the risk-averse agents, when the set of all Pareto optimal solutions, in the sense of viscosity, for the associated dynamic programming equations is dense in the given closed convex set. Moreover, for such a risk-averse decision problem, where results from the dynamic risk analysis are part of the information used in the risk-averse decision criteria, we briefly comment on the characteristics of acceptable risks vis-\'{a}-vis some uncertain future costs or outcomes.

Here, it is worth mentioning that some interesting studies on the dynamic risk measures, based on the conditional $g$-expecations, have been reported in the literature (e,g. see \cite{Pen04}, \cite{CorHMP02} and \cite{Ros06} for establishing connection between the risk measures and the generator of BSDE; and see also \cite{RSta10} for characterizing the generator of BSDE according to different risk measures). Moreover, such risk measures are widely used for evaluating the risk of uncertain future outcomes, and also assisting with stipulating minimum interventions for risk management (e.g., see \cite{ArtDEH99}, \cite{Pen04}, \cite{ElKPQ97}, \cite{FolS02}, \cite{DetS05} or \cite{CorHMP02} for related discussions). Recently, the authors in \cite{Rus10} and \cite{BefVP15} have provided interesting results on the risk-averse decision problem for Markov decision processes, in discrete-time setting, and, respectively, a hierarchical risk-averse framework for controlled-diffusion processes. Note that the rationale behind our framework follows in some sense the settings of these papers. However, to our knowledge, the problem of risk-aversion for controlled-diffusion processes has not been addressed in the context of multiple risk-averse agents argument, and it is important because it provides a mathematical framework that shows how a such framework can be systematically used to obtain consistently optimal risk-averse decisions.\footnote{In this paper, our intent is to provide a theoretical framework, rather than considering a specific numerical problem or application.}

The remainder of this paper is organized as follows. In Section~\ref{S2}, we present some preliminary results that are useful for our main results. In Section~\ref{S3}, using the basic remarks made in Section~\ref{S2}, we state the decision problem for the controlled-diffusion process with multiple risk-averse agents. In Section~\ref{S4}, we present our main results -- where we introduce a framework that requires a ``rational" cooperation among the risk-averse agents so as to achieve an overall optimal risk-averseness. Moreover, we establish the existence of optimal risk-averse solutions for the associated risk-averse dynamic programming equations. Finally, Section~\ref{S5} provides further remarks. 

\section{Preliminary results} \label{S2}
Let $\bigl(\Omega, \mathcal{F},\{\mathcal{F}_t \}_{t \ge 0}, \mathbb{P}\bigr)$ be a probability space, and let $\{B_t\}_{t \ge 0}$ be a $d$-dimensional standard Brownian motion, whose natural filtration, augmented by all $\mathbb{P}$-null sets, is denoted by $\{\mathcal{F}_t\}_{t \ge 0}$, so that it satisfies the {\it usual hypotheses} (e.g., see \cite{Pro90}). \,We consider the following controlled-diffusion process over a given finite-time horizon $T>0$
\begin{align}
d X_t^{u_{\cdot}} = m\bigl(t, X_t^{u_{\cdot}}, (\begin{array}{cccc} u_t^1, & u_t^2, & \cdots, & u_t^n \end{array})\bigr) dt +& \sigma\bigl(t, X_t^{u_{\cdot}}, (\begin{array}{cccc} u_t^1, & u_t^2, & \cdots, & u_t^n \end{array})\bigr)dB_t, \notag \\
& \quad\quad X_0^{u_{\cdot}}=x, \quad  0 \le t \le T,  \label{Eq1} 
\end{align}
where
\begin{itemize}
\item $X_{\cdot}^{u_{\cdot}}$ is an $\mathbb{R}^{d}$-valued controlled-diffusion process,
\item $u_{\cdot}^j$ is a $U^{j}$-valued measurable decision processes, which corresponds to the $j$th risk-averse agent (where $U^j$ is an open compact set in $\mathbb{R}^{m_j}$, with $j =1, 2, \ldots, n$); and, furthermore, $u_{\cdot} \triangleq (\begin{array}{cccc} u_{\cdot}^1, & u_{\cdot}^2, & \cdots, & u_{\cdot}^n \end{array})$ is an $n$-tuple of $\prod\nolimits_{i=1}^n U^{i}$-valued measurable decision processes such that for all $t > s$, $(B_t - B_s)$ is independent of $u_r$ for $r \le s$ (nonanticipativity condition) and
\begin{align*}
\mathbb{E} \int_{s}^{t} \vert u_{\tau}\vert^2 d\tau < \infty \quad \forall t \ge s,
\end{align*}
\item $m \colon [0, T] \times \mathbb{R}^d \times \prod_{i=1}^n U^{i} \rightarrow \mathbb{R}^{d}$ is uniformly Lipschitz, with bounded first derivative, and
\item $\sigma \colon [0, T] \times \mathbb{R}^{d} \times \prod\nolimits_{i=1}^n U^{i} \rightarrow \mathbb{R}^{d \times d}$ is Lipschitz with the least eigenvalue of $\sigma\,\sigma^T$ uniformly 
bounded away from zero for all $(x, u) \in \mathbb{R}^{d} \times \prod\nolimits_{i=1}^n U^{i}$ and $t \in [0, T]$, i.e., 
\begin{align*}
 \sigma(t, x, u)\,\sigma^T(t, x, u) \succeq \lambda I_{d \times d}, &\quad \forall (x, u) \in \mathbb{R}^{d} \times \prod\nolimits_{i=1}^n U^{i},\\
  & \quad \forall t \in [0, T],  
\end{align*}
for some $\lambda > 0$.
\end{itemize}

{\bf Notation}: Let us introduce the following spaces that will be useful later in the paper.
\begin{itemize}
\item $L^2\bigl(\Omega, \mathcal{F}_t, \mathbb{P}; \mathbb{R}^{d} \bigr)$ is the set of $\mathbb{R}^{d}$-valued $\mathcal{F}_t$-measurable random variables $\xi$ such that $\bigl\Vert \xi \bigr\Vert^2 = \mathbb{E}\bigl\{\bigl\vert \xi \bigr\vert^2  \bigr\}< \infty$;
\item $L^{\infty}\bigl(\Omega, \mathcal{F}_t, \mathbb{P}\bigr)$ is the set of $\mathbb{R}$-valued $\mathcal{F}_t$-measurable random variables $\xi$ such that $\bigl\Vert \xi \bigr\Vert = \essinf \bigl\vert \xi \bigr \vert < \infty$;
\item $\mathcal{S}^2\bigl(t, T; \mathbb{R}^{d} \bigr)$ is the set of $\mathbb{R}^{d}$-valued adapted processes $\bigl (\varphi_s\bigr)_{t \le s \le T}$ on $\Omega \times [t, T]$ such that $\bigl\Vert \varphi \bigr\Vert_{[t, T]}^2 = \mathbb{E}\bigl\{\sup_{t \le s \le T} \bigl\vert \varphi_s \bigr\vert^2  \bigr\}< \infty$;
\item $\mathcal{H}^2\bigl(t, T; \mathbb{R}^{d} \bigr)$ is the set of $\mathbb{R}^{d}$-valued progressively measurable processes $\bigl (\varphi_s\bigr)_{t \le s \le T}$ such that $\bigl\Vert \varphi \bigr\Vert_{[t, T]}^2 = \mathbb{E}\bigl\{ \int_t^T \bigl\vert \varphi_s \bigr \vert^2 ds  \bigr\}< \infty$.
\end{itemize}

On the same probability space $\bigl(\Omega, \mathcal{F},\{\mathcal{F}_t \}_{t \ge 0}, \mathbb{P}\bigr)$, we consider the following backward stochastic differential equation (BSDE) 
\begin{align}
- d Y_t = g\bigl(t, Y_t, Z_t\bigr) dt - Z_tdB_t, \quad Y_T=\xi, \label{Eq2}
\end{align}
where the terminal value $Y_T=\xi$ belongs to $L^2\bigl(\Omega, \mathcal{F}_T, \mathbb{P}; \mathbb{R}\bigr)$ and the generator functional $g \colon \Omega \times [0, T] \times \mathbb{R} \times \mathbb{R}^{d} \rightarrow \mathbb{R}$, with property that $\bigl(g\bigl(t, y, z\bigr)\bigr)_{0 \le t \le T}$ is progressively measurable for each $(y, z) \in \mathbb{R} \times \mathbb{R}^{d}$. We also assume that $g$ satisfies the following assumption.

\begin{assumption} \label{AS1}~\\\vspace{-5.0mm} 
\begin{enumerate} [{\rm (i)}]
\item $g$ is Lipschitz in $(y, z)$, i.e., there exists a constant $K > 0$ such that, $\mathbb{P}$-a.s., for any $t \in [0, T]$, $y_1, y_2 \in \mathbb{R}$ and $z_1, z_2 \in \mathbb{R}^d$ 
\begin{align*}
\bigl\vert g\bigl(t, y_1, z_1\bigr) - g\bigl(t, y_2, z_2\bigr) \bigr\vert \le K \bigl(\bigl\vert y_1 - y_2 \bigr\vert + \bigl\Vert z_1 - z_2 \bigr\Vert\bigr).
\end{align*}
\item $g\bigl(t, 0, 0\bigr) \in \mathcal{H}^2\bigl(t, T; \mathbb{R} \bigr)$.
\item $\mathbb{P}$-a.s., for all $t \in [0, T]$ and $y \in \mathbb{R}$, $g\bigl(t, y, 0\bigr) = 0$.
\end{enumerate}
\end{assumption}
Then, we state the following lemma, which is used to establish the existence of a unique adapted solution (e.g., see \cite{ParP90} for additional discussions).
\begin{lemma} \label{L1}
Suppose that Assumption~\ref{AS1} holds true. Then, for any $\xi \in L^2\bigl(\Omega, \mathcal{F}_T, \mathbb{P}; \mathbb{R}\bigr)$, the BSDE in \eqref{Eq2}, with terminal condition $Y_T=\xi$, i.e.,
\begin{align}
 Y_t = \xi + \int_t^T g\bigl(s, Y_s, Z_s\bigr) ds - \int_t^TZ_s dB_s, \quad 0 \le t \le T  \label{Eq3}
\end{align}
has a unique adapted solution
\begin{align}
 \bigl(Y_t^{T,g,\xi}, Z_t^{T,g,\xi}\bigr)_{0 \le t \le T} \in \mathcal{S}^2\bigl(0, T; \mathbb{R} \bigr) \times  \mathcal{H}^2\bigl(0, T; \mathbb{R}^{d} \bigr). \label{Eq4}
\end{align}
\end{lemma}

Moreover, we recall the following comparison theorem, which is restricted to one-dimensional BSDEs (e.g., see \cite{ParT99}).
\begin{theorem} \label{T1}
Given two generators $g_1$ and $g_2$ satisfying Assumption~\ref{AS1} and two terminal conditions $\xi_1, \xi_2 \in L^2\bigl(\Omega, \mathcal{F}_T, \mathbb{P}; \mathbb{R}\bigr)$. Let $\bigl(Y_t^1, Z_t^1\bigr)$ and $\bigl(Y_t^2, Z_t^2\bigr)$ be the solution pairs corresponding to $\bigl(\xi_1, g_1\bigr)$ and $\bigl(\xi_2, g_2\bigr)$, respectively. Then, we have 
\begin{enumerate} [{\rm (i)}]
\item Monotonicity: If $\xi_1 > \xi_2$ and $g_1 > g_2$, $\mathbb{P}$-a.s., then $Y_t^1 > Y_t^2$, $\mathbb{P}$-a.s., for all $t \in [0, T]$;
\item Strictly Monotonicity: In addition to ($i$) above, if we assume that $\mathbb{P}\bigl(\xi_1 > \xi_2\bigr) > 0$, then $\mathbb{P}\bigl(Y_t^1 > Y_t^2\bigr) > 0$, for all $t \in [0, T]$.
\end{enumerate}
\end{theorem}
In the following, we give the definition for a dynamic risk measure that is associated with the generator of BSDE in \eqref{Eq2}.

\begin{definition} \label{Df1} 
For any $\xi \in L^2\bigl(\Omega, \mathcal{F}_T, \mathbb{P}; \mathbb{R}\bigr)$, let $\bigl(Y_t^{T,g,\xi}, Z_t^{T,g,\xi}\bigr)_{0 \le t \le T} \in  \mathcal{S}^2\bigl(0, T; \mathbb{R} \bigr) \times  \mathcal{H}^2\bigl(0, T; \mathbb{R}^{d} \bigr)$ be the unique solution for the BSDE in \eqref{Eq2} with terminal condition $Y_T=\xi$. Then, we define the dynamic risk measure $\rho_{t,T}^g$ of $\xi$ by\footnote{Here, we remark that, for any $t \in [0,T]$, the conditional $g$-expectation (denoted by $\mathcal{E}_g\bigl[\xi \vert \mathcal{F}_t\bigr]$) is also defined by
\begin{align*}
\mathcal{E}_g\bigl[\xi \vert \mathcal{F}_t\bigr] \triangleq Y_t^{T,g,\xi}.
\end{align*}} 
\begin{align}
\rho_{t,T}^g \bigl[\xi \bigr] \triangleq Y_t^{T,g,\xi}.  \label{Eq5}
\end{align}
\end{definition}

Moreover, if the generator functional $g$ satisfies Assumption~\ref{AS1}, then a family of time-consistent dynamic risk measures $\bigl\{\rho_{t,T}^g\bigr\}_{t \in [0,T]}$ has the following properties (see \cite{Pen04} for additional discussions).
 
\begin{property} ~ \\ \vspace{-5.0mm}
\begin{enumerate} [(i)]
\item {\it Convexity}: If $g$ is convex for every fixed $(t, \omega) \in [0, T] \times \Omega$, then for all $\xi_1, \xi_2 \in L^2\bigl(\Omega, \mathcal{F}_T, \mathbb{P}; \mathbb{R}\bigr)$ and for all $\lambda \in L^{\infty}\bigl(\Omega, \mathcal{F}_t, \mathbb{P}; \mathbb{R}\bigr)$ such that $0 \le \lambda \le 1$
\begin{align*}
\rho_{t,T}^g\bigl[\lambda \xi_1 + (1-\lambda)\xi_2 \bigr] \le \lambda \rho_{t,T}^g\bigl[\xi_1\bigr] + (1- \lambda) \rho_{t,T}^g\bigl[\xi_1\bigr];
\end{align*}
\item {\it Monotonicity}:  For $\xi_1, \xi_2 \in L^2\bigl(\Omega, \mathcal{F}_T, \mathbb{P}; \mathbb{R}\bigr)$ such that $\xi_1 > \xi_2$ $\mathbb{P}$-a.s., then 
\begin{align*}
\rho_{t,T}^g\bigl[\xi_1\bigr] > \rho_{t,T}^g\bigl[\xi_2\bigr], \quad \mathbb{P}{\text-a.s.};
\end{align*}
\item {\it Trans-invariance}: For all $\xi \in L^2\bigl(\Omega, \mathcal{F}_T, \mathbb{P}; \mathbb{R}\bigr)$ and $\nu \in L^2\bigl(\Omega, \mathcal{F}_t, \mathbb{P}; \mathbb{R}\bigr)$
\begin{align*}
\rho_{t,T}^g\bigl[\xi + \nu\bigr] = \rho_{t,T}^g\bigl[\xi\bigr] + \nu;
\end{align*}
\item {\it Positive-homogeneity}: For all $\xi \in L^2\bigl(\Omega, \mathcal{F}_T, \mathbb{P}; \mathbb{R}\bigr)$ and for all $\lambda \in L^{\infty}\bigl(\Omega, \mathcal{F}_t, \mathbb{P}; \mathbb{R}\bigr)$ such that $\lambda > 0$
\begin{align*}
\rho_{t,T}^g\bigl[\lambda \xi \bigr] = \lambda \rho_{t,T}^g\bigl[\xi\bigr];
\end{align*}
\item {\it Normalization}: $\rho_{t,T}^g\bigl[0\bigr] = 0$ for $t \in [0, T]$.
\end{enumerate}
\end{property}

\begin{remark} \label{R1.1}
In Section~\ref{S3}, using the basic remarks made in this section, we introduce multi-structure dynamic risk measures (that satisfy the above properties {\rm (i)}--{\rm (v)}) induced from conditional $g$-expectations, where the latter are associated with the generator functionals of certain BSDEs that implicitly take into account the risk-cost functionals of the risk-averse agents.
\end{remark}

In this paper, we consider a risk-averse decision problem for the above controlled-diffusion process, in which the decision makers (i.e., the $n$ {\it risk-averse agents} with differing risk-averse related responsibilities and information) choose their risk-averse decisions from progressively measurable strategy sets. That is, the $j$th-agent's decision $u_{\cdot}^j$ is a $U^j$-valued measurable control process from
\begin{align}
\mathcal{U}_{[0,T]}^j \triangleq  \Bigl\{u^j \colon [0,T] \times & \Omega \rightarrow U^j \,\Bigl\vert \, u^j \,\, \text{is an} \,\, \bigl\{\mathcal{F}_t\bigr\}_{t\ge 0}\text{- adapted} \notag \\
& \quad \text{and}\,\, \mathbb{E} \int_{0}^{T} \vert u_t^j\vert^2 dt < \infty \Bigr\}, \quad j=1,2, \dots, n. \label{Eq6} 
\end{align}
Moreover, we suppose that the risk-averse agents are ``rational" (in the sense of making consistent decisions that minimize their individual accumulated risk-costs) with a certain $n$-tuple of measurable decision processes $\hat{u} =  (\begin{array}{cccc} \hat{u}_{\cdot}^1, & \hat{u}_{\cdot}^2, & \cdots, & \hat{u}_{\cdot}^n \end{array}) \in \bigotimes\nolimits_{i=1}^n \mathcal{U}_{[0,T]}^i$. Further, we consider the following cost functionals providing information about the accumulated risk-costs on the time interval $[0, T]$ w.r.t. each of the risk-averse agents, i.e., 
\begin{align}
 \xi_{0,T}^j(u^{\neg j}) = \int_0^T c_j\bigl(t, X_t^{u_{\cdot}^{\neg j}}, u_t^j\bigr) dt + \Psi_j(X_T^{u_{\cdot}^{\neg j}}), \quad j=1,2, \dots, n, \label{Eq7} 
\end{align}
where 
\begin{align*}
u_{\cdot}^{\neg j} \triangleq (\begin{array}{ccccccc} \hat{u}_{\cdot}^1, & \cdots, & \hat{u}_{\cdot}^{j-1}, & u_{\cdot}^j, & \hat{u}_{\cdot}^{j+1}, & \cdots, & \hat{u}_{\cdot}^n \end{array}) \in \bigotimes\nolimits_{i=1}^n \mathcal{U}_{[0,T]}^i, 
\end{align*}
with $c_j \colon [0,T] \times \mathbb{R}^d \times U^j \rightarrow \mathbb{R}$ and $\Psi_j\colon \mathbb{R}^d \rightarrow \mathbb{R}$ are measurable functions. Here, we remark that the corresponding solution $X_t^{u_{\cdot}^{\neg j}}$, for $j \in \{1,2, \dots, n\}$, in equation~\eqref{Eq1} depends on the $n$-tuple admissible risk-averse decisions $u_{\cdot}^{\neg j} \in \bigotimes\nolimits_{i=1}^n \mathcal{U}_{[0,T]}^i$ and it also depends on the initial condition $X_0^{u_{\cdot}^{\neg j}}=x$.  As a result of this, for any time-interval $[t, T]$, with $t \in [0, T]$, the accumulated risk-costs $ \xi_{t,T}^j$, for $j=1,2, \dots, n$, depend on the risk-averse decisions $u_{\cdot}^{\neg j} \in \bigotimes\nolimits_{i=1}^n \mathcal{U}_{[t,T]}^i$.\footnote{Here, we use the notation $u^{\neg j}$ to emphasize the dependence on $u_{\cdot}^j \in \mathcal{U}_{[t,T]}^j$, where $\mathcal{U}_{[t,T]}^j$, for any $t \in [0, T]$, denotes the sets of $U^j$-valued $\bigl\{\mathcal{F}_s^t\bigr\}_{s \ge t}$-adapted processes (see Definition~\ref{Df2}).}  Moreover, we also assume that $f$, $\sigma$, $c_j$ and $\Psi_j$, for $p \ge 1$, satisfy the following growth conditions  
\begin{align}
\bigl\vert m\bigl(t, x, u) \bigr\vert &+ \bigl\vert \sigma\bigl(t, x, u \bigr) \bigr\vert + \bigl\vert c_j\bigl(t, x, u\bigr) \bigr\vert + \bigl\vert \Psi_j\bigl(x\bigr) \bigr\vert \notag \\
 & \quad \le K \bigl(1 + \bigl\vert x \bigr\vert^p + \bigl\vert u \bigr\vert \bigr), \quad \forall j \in \{1,2, \dots, n\}, \label{Eq8}
\end{align}
for all $\bigl(t, x, u \bigr) \in [0,T] \times \mathbb{R}^{d} \times \prod_{i=1}^n U^i$ and for some constant $K > 0$.

\section{Risk-averse decision problem formulation} \label{S3}
In order to make our problem formulation more precise, for any $(t, x) \in [0, T] \times \mathbb{R}^d$, we consider the following forward-SDE with an initial condition $X_t^{t,x;u_{\cdot}^{\neg j}} = x$, for $j \in \{1,2, \ldots, n\}$,
\begin{align}
d X_s^{t,x;u_{\cdot}^{\neg j}} &= m\bigl(s, X_s^{t,x;u_{\cdot}^{\neg j}}, (\begin{array}{ccccccc} \hat{u}_s^1, & \cdots, & \hat{u}_s^{j-1}, & u_s^{j}, & \hat{u}_s^{j+1},& \cdots, & \hat{u}_s^n \end{array})\bigr) dt \notag \\
& \quad \quad + \sigma\bigl(s, X_s^{t,x;u_{\cdot}^{\neg j}}, (\begin{array}{ccccccc} \hat{u}_s^1, & \cdots, & \hat{u}_s^{j-1}, & u_s^{j}, & \hat{u}_s^{j+1},& \cdots, & \hat{u}_s^n \end{array})\bigr)dB_t, \notag \\
 & \quad \quad \quad \quad \quad \quad X_t^{t,x;u_{\cdot}^{\neg j}} = x, \quad t \le s \le T,  \label{Eq9}
\end{align}
where $u_{\cdot}^{\neg j} = (\begin{array}{ccccccc} \hat{u}_s^1, & \cdots, & \hat{u}_s^{j-1}, & u_s^{j}, & \hat{u}_s^{j+1},& \cdots, & \hat{u}_s^n \end{array})$ is an $n$-tuple of $\prod\nolimits_{i=1}^n U^j$-valued measurable decision processes.

Let $\bigl\{\xi_j^{Target}\bigr\}_{j=1}^n$ be a set of real-valued random variables from $L^2(\Omega, \mathcal{F}_T, \mathbb{P}; \mathbb{R})$ and we further suppose that the data $\xi_j^{Target}$ take the following forms
\begin{align}
 \xi_j^{Target} = \Psi_j(X_T^{t,x;u_{\cdot}^{\neg j}}),\quad  j=1,2, \ldots, n, \quad \mathbb{P}-a.s.  \label{Eq10}
 \end{align}
Moreover, we introduce the following risk-value functions
\begin{align}
 V_j^{u^j}\bigl(t, x\bigr) = \rho_{t, T}^{g_j} \bigl[\xi_{t,T}^j\bigl(u^{\neg j}\bigr)\bigr], \quad j=1,2, \ldots, n,   \label {Eq11}
\end{align}
where
\begin{align}
\xi_{t,T}^j \bigl(u^{\neg j}\bigr) = \int_t^T c_j\bigl(s, X_s^{t,x;u_{\cdot}^{\neg j}}, u_s^j\bigr) ds + \Psi_j(X_T^{t,x;u_{\cdot}^{\neg j}}). \label{Eq12}
\end{align}
Then, taking into account equation~\eqref{Eq10} (and with the Markovian framework), we can express the above risk-value functions using standard-BSDEs as follows
\begin{align}
 V_j^{u^j}\bigl(t, x\bigr) &\triangleq Y_s^{j,t,x;u_{\cdot}^{\neg j}} \notag \\
                             &= \Psi_j(X_T^{t,x;u_{\cdot}^{\neg j}}) + \int_t^T g_j\bigl(s, X_s^{t,x;u_{\cdot}^{\neg j}} Y_s^{j, t,x;u_{\cdot}^{\neg j}}, Z_s^{j,t,x;u_{\cdot}^{\neg j}}\bigr)ds \notag \\
                                                                      & \quad \quad \quad - \int_t^T Z_s^{j,t,x; u_{\cdot}^{\neg j}}dB_s, \quad j=1,2, \dots, n, \label{Eq13} 
\end{align}
where 
\begin{align*}
& g_j\bigl(t, X_s^{t,x;u_{\cdot}^{\neg j}}, Y_s^{j, t,x;u_{\cdot}^{\neg j}}, Z_s^{j, t,x;w}\bigr) \notag \\
&  \quad \quad \quad \quad \quad \quad= c_j\bigl(s, X_s^{t,x;u_{\cdot}^{\neg j}}, u_s^j\bigr) + g\bigl(s,Y_s^{j, t,x;u_{\cdot}^{\neg j}}, Z_s^{j,t,x;u_{\cdot}^{\neg j}}\bigr).
\end{align*}
and further noting the conditions in \eqref{Eq8}, then the pairs $\bigl(Y_s^{j,t,x;u_{\cdot}^{\neg j}}, Z_s^{j,t,x;u_{\cdot}^{\neg j}}\bigr)_{t \le s \le T}$ are adapted solutions on $[t, T] \times \Omega$ and belong to $\mathcal{S}^2\bigl(t, T; \mathbb{R} \bigr) \times  \mathcal{H}^2\bigl(t, T; \mathbb{R}^{d} \bigr)$. Equivalently, we can also rewrite \eqref{Eq13} as a family of BSDEs on the probability space $\bigl(\Omega, \mathcal{F}, \mathbb{P}, \{\mathcal{F}_t \}_{t \ge 0})$, i.e., for $s \in [t, T]$,
\begin{align}
-d Y_s^{j,t,x;u_{\cdot}^{\neg j}} = g_j\bigl(s, X_s^{t,x;u_{\cdot}^{\neg j}} Y_s^{j, t,x;u_{\cdot}^{\neg j}}, Z_s^{j,t,x;u_{\cdot}^{\neg j}}\bigr) ds - Z_s^{j,t,x; u_{\cdot}^{\neg j}}dB_s, \notag \\
    Y_T^{j,t,x;u_{\cdot}^{\neg j}} =  \Psi_j(X_T^{t,x;u_{\cdot}^{\neg j}}), \quad j=1,2, \dots, n. \label{Eq14}
\end{align}
In the following, we denote the solutions $\bigl(\begin{array}{cccc} Y_s^{1,t,x;u_{\cdot}^{\neg 1}}, & Y_s^{2,t,x;u_{\cdot}^{\neg 2}}, & \cdots, & Y_s^{n,t,x;u_{\cdot}^{\neg n}} \end{array}\bigr) \in \bigotimes\nolimits_{i=1}^n \mathcal{S}^2\bigl(t, T; \mathbb{R} \bigr)$ and $\bigl(\begin{array}{cccc} Z_s^{1,t,x;u_{\cdot}^{\neg 1}}, & Z_s^{2,t,x;u_{\cdot}^{\neg 2}}, & \cdots, & Z_s^{n,t,x;u_{\cdot}^{\neg n}} \end{array}\bigr) \in  \bigotimes\nolimits_{i=1}^n \mathcal{H}^2\bigl(t, T; \mathbb{R}^d \bigr)$ by bold letters $\mathbf{Y}_s^{t,x;u}$ and $\mathbf{Z}_s^{t,x;u}$, respectively, for any $t \in [0, T]$ and for $s \in [t, T]$.  Similarly, the family of BSDEs in \eqref{Eq14} can be rewritten as a multi-dimensional BSDE as follows
\begin{align}
-d \mathbf{Y}_s^{t,x;u} = \mathbf{G}\bigl(s, X_s^{t,x;u} \mathbf{Y}_s^{t,x;u}, \mathbf{Z}_s^{t,x;u}\bigr) ds &- \mathbf{Z}_s^{t,x; u} dB_s, \notag \\
    s \in [t, T], \quad &\mathbf{Y}_T^{t,x;u} =  \mathbf{\Psi}(X_T^{t,x;u}), \label{Eq15}
\end{align}
where
\begin{align*}
&\mathbf{G}\bigl(s, X_s^{t,x;u} \mathbf{Y}_s^{t,x;u}, \mathbf{Z}_s^{t,x;u}\bigr) \\
&\quad = \blockdiag \Bigl\{ g_1\bigl(s, X_s^{t,x;u_{\cdot}^{\neg 1}} Y_s^{1, t,x;u_{\cdot}^{\neg 1}}, Z_s^{1,t,x;u_{\cdot}^{\neg 1}}\bigr),\\
& \quad \quad g_2\bigl(s, X_s^{t,x;u_{\cdot}^{\neg 2}} Y_s^{2, t,x;u_{\cdot}^{\neg 2}}, Z_s^{2,t,x;u_{\cdot}^{\neg 2}}), \cdots, g_n\bigl(s, X_s^{t,x;u_{\cdot}^{\neg n}} Y_s^{n, t,x;u_{\cdot}^{\neg n}}, Z_s^{n,t,x;u_{\cdot}^{\neg n}}\bigr) \bigr)\Bigr\}
\end{align*}
and
\begin{align*}
\mathbf{\Psi}(X_T^{t,x;u}) = \bigl(\begin{array}{cccc} \Psi_1(X_T^{t,x;u_{\cdot}^{\neg 1}}), & \Psi_2(X_T^{t,x;u_{\cdot}^{\neg 2}}), & \cdots, & \Psi_n(X_T^{t,x;u_{\cdot}^{\neg n}}) \end{array}\bigr).
\end{align*}

Let $K$ be a closed convex set in $\mathbb{R}^n$, then we recall the notion of viability property for the BSDE in \eqref{Eq15} (cf. equations~\eqref{Eq13} and \eqref{Eq14}).
\begin{definition} \label{Df2}
Let  $\hat{u}_{\cdot} = (\begin{array}{cccc} \hat{u}_{\cdot}^1, & \hat{u}_{\cdot}^2, & \cdots, & \hat{u}_{\cdot}^n \end{array}) \in \bigotimes\nolimits_{i=1}^n \mathcal{U}_{[t,T]}^i$ be an $n$-tuple of ``rational" preferable decisions for the risk-averse agents. Then, for a nonempty closed convex set $K \subset \mathbb{R}^n$ and for $u_{\cdot}^j \in \mathcal{U}_{[0,T]}^j$, with $j=1,2, \ldots, n$
\begin{enumerate} [(a)]
\item A stochastic process $\bigl\{\mathbf{Y}_t^{0,x;u}, \,\,t \in [0, T] \bigr\}$ is viable in $K$ if and only if for $\mathbb{P}$-{\it almost}\, $\omega \in \Omega$ 
\begin{align}
\mathbf{Y}_t^{0,x;u}(\omega) \in K, \quad \forall t \in [0, T]. \label{Eq16}
\end{align}
\item The closed convex set $K$ enjoys the Backward Stochastic Viability Property (BSVP) for the equation in \eqref{Eq15} if and only if for all $\tau \in [0, T]$, with equation~\eqref{Eq10}, i.e., 
\begin{align}
\forall \,\mathbf{\Xi}^{Target} = \bigl(\begin{array}{cccc} \xi_1^{Target}, & \xi_2^{Target}, & \cdots, & \xi_n^{Target}\end{array}\bigr) \in L^2\bigl(\Omega, \mathcal{F}_{\tau}, \mathbb{P}; \mathbb{R}^{n} \bigr), \label{Eq17}
\end{align}
there exists a solution pair $\bigl(\mathbf{Y}_{\cdot}^{0,x;u}, \mathbf{Z}_{\cdot}^{0,x;u} \bigr)$ to the BSDE in \eqref{Eq15} over the time interval $[0, \tau]$,
\begin{align}
\mathbf{Y}_s^{0,x;u} &= \mathbf{\Xi}^{Target} + \int_s^{\tau} \mathbf{G}\bigl(r, X_r^{0,x;u} \mathbf{Y}_r^{0,x;u}, \mathbf{Z}_r^{0,x;u}\bigr) dr -  \int_s^{\tau} \mathbf{Z}_r^{0,x; u} dB_r, \label{Eq18}
\end{align}
with 
\begin{align*}
\bigl(\mathbf{Y}_{\cdot}^{0,x;u}, \mathbf{Z}_{\cdot}^{0,x;u} \bigr) \in \bigotimes\nolimits_{i=1}^n  \mathcal{S}^2\bigl(0, \tau; \mathbb{R} \bigr) \times \bigotimes\nolimits_{i=1}^n \mathcal{H}^2\bigl(0, \tau; \mathbb{R}^{d} \bigr),
\end{align*}
such that $\bigl\{\mathbf{Y}_s^{0,x;u}, \,\,s \in [0, \tau] \bigr\}$ is viable in $K$.
\end{enumerate}
\end{definition}

For the above given closed convex set $K$, let us define the projection of a point $a$ onto $K$ as follow
\begin{align}
\Pi_{K} (a) = \Bigl\{b \in K \, \bigl\vert \, \vert a - b \vert = \min _{c \in K} \vert a - c \vert = d_{K}(a) \Bigr\}. \label{Eq19}
\end{align}
Notice that, since $K$ is convex, from the Motzkin's theorem, $\Pi_{K}$ is single-valued. Further, we recall that $d_{K}^2(\cdot)$ is convex; and thus, due to Alexandrov's theorem \cite{Ale39}, $d_{K}^2(\cdot)$ is almost everywhere twice differentiable.

Assume that there exists an $n$-tuple of ``rational" decisions $\hat{u} =  (\begin{array}{cccc} \hat{u}_{\cdot}^1, & \hat{u}_{\cdot}^2, & \cdots, & \hat{u}_{\cdot}^n \end{array}) \in \bigotimes\nolimits_{i=1}^n \mathcal{U}_{[0,T]}^i$ which is preferable by all risk-averse decision-making agents. Moreover, on the space $C_b^{1,2}([t, T] \times \mathbb{R}^d); \mathbb{R}^n)$, for any $(t,x) \in [0, T] \times \mathbb{R}^d$, we consider the following system of semilinear parabolic partial differential equations (PDEs)
\begin{eqnarray}
\left.\begin{array}{r}
 \dfrac{\partial \varphi_j(t, x)}{\partial t}  + \inf_{u^j \in U^j} \Bigl\{\mathcal{L}_{t}^{u^{\neg j}} \varphi_j(t, x) ~~~ \hspace{1.5in} \quad \quad \quad  \\
  + g_j\bigl(t, \varphi(t, x), D_x \varphi_j(t, x) \cdot \sigma(t, x, u^{\neg j})\bigr)\Bigr\} = 0  \\
   j = 1, 2, \ldots, n 
\end{array}\right\}  \label{Eq20}
\end{eqnarray}
with the following boundary condition
\begin{align}
\varphi(T, x) &= \mathbf{\Psi}(x), \notag \\
                    &\equiv \bigl(\begin{array}{cccc} \Psi_1(x), & \Psi_2(x), & \cdots, & \Psi_n(x) \end{array}\bigr), \quad x \in \mathbb{R}^d,  \label{Eq21}
\end{align}
where, for any $\phi(x) \in C_0^{\infty}(\mathbb{R}^d)$, the second-order linear operators $\mathcal{L}_{t}^{u^{\neg j}}$ are given by
\begin{align}
 \mathcal{L}_{t}^{u^{\neg j}} \phi(x) = \dfrac{1}{2} \operatorname{tr} \Bigl\{a(t, x, u^{\neg j}) D_{x}^2 \phi(x)\Bigr\} + m(t, x, u^{\neg j}) D_{x} \phi(x), \,\, \notag \\
                                            t \in [0, T], \quad j =1, 2, \ldots n, \label{Eq22}
\end{align}
with $a(t, x, u^{\neg j}) = \sigma(t, x, u^{\neg j}) \sigma^T(t, x, u^{\neg j})$, $D_{x}$ and $D_{x}^2$, (with $D_{x}^2 = \bigl({\partial^2 }/{\partial x_k \partial x_l} \bigr)$) are the gradient and the Hessian (w.r.t. the variable $x$), respectively.

\begin{remark} \label{R3.1}
Here, we remark that the above system of equations in \eqref{Eq20} together with \eqref{Eq21}, is associated with the decision problem for the risk-averse agents, restricted to $\Sigma_{[t,T]}$ (see Definition~\ref{Df2} below). Moreover, such a system of equations represents a generalized family of HJB equation with additional terms $g_j$ for $j = 1,2, \ldots, n$. Note that the problem of FBSDEs (cf. equations~\eqref{Eq9} and \eqref{Eq15} or \eqref{Eq14}) and the solvability of the related system of semilinear parabolic PDEs have been well studied in literature (e.g., see \cite{Ant93}, \cite{HUP95}, \cite{LIW14}, \cite{ParT99}, \cite{Pen91} and \cite{Pen92}).
\end{remark}

Next, we recall the definition of viscosity solutions for \eqref{Eq20} along with \eqref{Eq21} (e.g., see \cite{CraIL92}, \cite{FleS06} or \cite{Kry08} for additional discussions on the notion of viscosity solutions).

\begin{definition} \label{Df3}
The function $\varphi \colon [0, T] \times \mathbb{R}^d \rightarrow \mathbb{R}^n$ is a viscosity solution for \eqref{Eq20} together with the boundary condition in \eqref{Eq21}, if the following conditions hold
\begin{enumerate} [(i)]
\item for every $\psi \in C_b^{1,2}([0, T], \times \mathbb{R}^d; \mathbb{R}^n)$ such that $\psi \ge \varphi$ on $[0, T] \times \mathbb{R}^d$,
\begin{align}
\sup_{(t,x)} \bigl\{\varphi(t,x) - \psi(t,x) \bigr\} = 0, \label{Eq23}
\end{align}
and for $(t_{0},x_{0}) \in  [0, T] \times \mathbb{R}^d$ such that  $\psi(t_{0}, x_{0})=\varphi(t_{0}, x_{0})$ (i.e., a local maximum at $(t_{0},x_{0})$), then we have
\begin{align}
 & \dfrac{\partial \psi_j(t_{0},x_{0})}{\partial t}  + \inf_{u^j \in U^j} \Bigl\{\mathcal{L}_{t}^{u^{\neg j}} \psi_j(t_{0},x_{0}) \notag \\
 & \quad \quad + g_j\bigl(t_{0}, x_{0}, \psi(t_{0},x_{0}), D_x \psi_j(t_{0},x_{0}) \cdot \sigma(t_{0},x_{0}, u^{\neg j})\bigr)\Bigr\} \ge 0 \label{Eq24}
\end{align}
\item for every $\psi \in C_b^{1,2}([0, T], \times \mathbb{R}^d; \mathbb{R}^n)$ such that $\psi \le \varphi$ on $[0, T] \times \mathbb{R}^d$,
\begin{align}
\inf_{(t,x)} \bigl\{ \varphi(t,x) - \psi(t,x) \bigr\} = 0,  \label{Eq3.11}
\end{align}
and for $(t_{0},x_{0}) \in  [0, T] \times \mathbb{R}^d$ such that  $\psi(t_{0}, x_{0})=\varphi(t_{0}, x_{0})$ (i.e., a local minimum at $(t_{0},x_{0})$), then we have
\begin{align}
 & \dfrac{\partial \psi_j(t_{0},x_{0})}{\partial t}  + \inf_{u^j \in U^j} \Bigl\{ \mathcal{L}_{t}^{u^{\neg j}} \psi_j(t_{0},x_{0}) \notag \\
 & \quad \quad+ g_j\bigl(t_{0},x_{0}, \psi(t_{0},x_{0}), D_x \psi_j(t_{0},x_{0}) \cdot \sigma(t_{0},x_{0}, u^{\neg j})\bigr)\Bigr\} \le 0, \label{Eq25}
\end{align}
\end{enumerate}
for $j = 1,2, \ldots, n$.
\end{definition}

Next, let us define the viability property for the system of semilinear parabolic PDEs in \eqref{Eq20} as follow.
\begin{definition} \label{Df4}
The system of semilinear parabolic PDEs in \eqref{Eq20} enjoys the viability property w.r.t. the closed convex set $K$ if and only if, for any $\mathbf{\Psi} \in C_p(\mathbb{R}^d; \mathbb{R}^n)$ taking values in $K$, the viscosity solution to \eqref{Eq20} satisfies
\begin{align}
\forall (t, x) \in [0, T] \times \mathbb{R}^d, \quad \varphi(t,x) \in K. \label{Eq27}
\end{align}
\end{definition}

Later in Section~\ref{S4}, assuming the Markovian framework, we provide additional results that establish a connection between the viability property of the BSDE in \eqref{Eq15}, w.r.t. the closed convex set $K$, and the solutions, in the sense viscosity, for the system of semilinear parabolic PDEs in \eqref{Eq20}.

In what follows, we introduce a framework that requires a ``rational" cooperation among the risk-averse agents so as to achieve an overall risk-averseness (in the sense of Pareto optimality). For example, for any $t \in [0, T]$, let us assume that 
\begin{align*}
\hat{u}_{\cdot} = (\begin{array}{cccc} \hat{u}_{\cdot}^1, & \hat{u}_{\cdot}^2, & \cdots, & \hat{u}_{\cdot}^n \end{array}) \in \bigotimes\nolimits_{i=1}^n \mathcal{U}_{[t,T]}^i
\end{align*}
is an $n$-tuple of ``rational" preferable decisions for the risk-averse agents, then the problem of finding an optimal risk-averse decision for the $j$th-agent, where $j \in \{1,2, \ldots, n\}$, that minimizes the $j$th-accumulated risk-cost functional, is equivalent to finding an optimal solution for
\begin{align}
\inf_{u_{\cdot}^j \in \mathcal{U}_{[t,T]}^j}  J_j\bigr[\bigl(u^{\neg j})\bigl], \label{Eq28}
\end{align}
where
\begin{align}
J_j\bigr[\bigl(u^{\neg j}\bigr)\bigl] = \rho_{t, T}^{g_j} \bigl[\xi_{t,T}^j\bigl(u^{\neg j}\bigr)\bigr], \label{Eq29}
\end{align}
with $u_{\cdot}^{\neg j} = (\begin{array}{ccccccc} \hat{u}_{\cdot}^1, &\ldots, & \hat{u}_{\cdot}^{j-1}, & u_{\cdot}^j, & \hat{u}_{\cdot}^{j+1}, & \cdots, & \hat{u}_{\cdot}^n \end{array}) \in \bigotimes\nolimits_{i=1}^n \mathcal{U}_{[t,T]}^i$.

\begin{remark} \label{R2.2}
Here, we remark that the generator functionals $g_j$, for $j=1,2, \ldots, n$, contain a common term $g$ that acts on different processes (see equations~\eqref{Eq13} and \eqref{Eq14}). Moreover, due to differing risk-cost functionals w.r.t. each of the agents, we also observe that $\bigl\{\rho_{t, T}^{g_j} \bigl[\,\cdot\,]\bigr\}_{j=1}^n$, for $t \in [0, T]$, in equation~\eqref{Eq29} provide multi-structure, time-consistent, dynamic risk measures vis-\'{a}-vis some uncertain future outcomes specified by a set of random variables from $L^2(\Omega, \mathcal{F}_T, \mathbb{P}; \mathbb{R})$.
\end{remark}

Note that, for any  given $u_{\cdot}^j \in \mathcal{U}_{[t,T]}^j$, if the forward-backward stochastic differential equations (FBSDEs) in \eqref{Eq9} and \eqref{Eq15} (cf. equations~\eqref{Eq13} and \eqref{Eq14}) admit unique solutions and, further, $\mathbf{Y}_s^{t,x;u} (\omega) \in K$, for $\mathbb{P}$-{\it almost} $\omega \in \Omega$ and for all $s \in [t, T]$ and for $t \in [0,T]$. Then, any ``rational" preferable decisions for the $j$th-agent satisfy the following
\begin{align}
  \hat{u}_{\cdot}^j \in \Bigl\{\tilde{u}_{\cdot} \in \mathcal{U}_{[t,T]}^j \,\Bigl\vert\, \rho_{t, T}^{g_j} \bigl[\xi_{t,T}^j\bigl(\tilde{u}^{\neg j}\bigr)\bigl] \le \rho_{t, T}^{g_j} \bigl[\xi_{t,T}^j\bigl(u^{\neg j}\bigr)\bigl], \quad \quad & \notag \\
   \forall (\begin{array}{ccccccc} \hat{u}_{\cdot}^1, &\ldots, & \hat{u}_{\cdot}^{j-1}, & \hat{u}_{\cdot}^{j+1}, & \cdots, & \hat{u}_{\cdot}^n \end{array}) \in \bigotimes\nolimits_{i \neq j} \mathcal{U}_{[t,T]}^i, \quad & \notag \\
   \quad  \forall j \in \{1,2, \ldots, n\},  \quad \mathbb{P}-a.s. \Bigr\}, & \label{Eq30}
\end{align}
where $\tilde{u}_{\cdot}^{\neg j} = (\begin{array}{ccccccc} \hat{u}_{\cdot}^1, &\ldots, & \hat{u}_{\cdot}^{j-1}, & \tilde{u}_{\cdot}^j, & \hat{u}_{\cdot}^{j+1}, & \cdots, & \hat{u}_{\cdot}^n \end{array}) \in \bigotimes\nolimits_{i=1}^n \mathcal{U}_{[t,T]}^i$.\footnote{In the paper, we assume that the set on the right-hand side of \eqref{Eq30} is nonempty.}

Next, we introduce the following definition for an admissible risk-averse decision system $\Sigma_{[t, T]}$, with multi-structure dynamic risk measures, which provides a logical construct for our main results (e.g., see also \cite{LIW14}).
\begin{definition} \label{Df5}
For a given finite-time horizon $T>0$, we call $\Sigma_{[t, T]}$ an admissible risk-averse decision system, if it satisfies the following conditions:
\begin{itemize}
\item $\bigl(\Omega, \mathcal{F},\{\mathcal{F}_t \}_{t \ge 0}, \mathbb{P}\bigr)$ is a complete probability space;
\item $\bigl\{B_s\bigr\}_{s \ge t}$ is a $d$-dimensional standard Brownian motion defined on $\bigl(\Omega, \mathcal{F}, \mathbb{P}\bigr)$ over $[t, T]$ and $\mathcal{F}^t \triangleq \bigl\{\mathcal{F}_s^t\bigr\}_{s \in [t, T]}$, where $\mathcal{F}_s^t = \sigma\bigl\{\bigl(B_s; \,t \le s \le T \bigr)\bigr\}$ is augmented by all $\mathbb{P}$-null sets in $\mathcal{F}$;
\item $u_{\cdot}^j \colon \Omega \times [s, T]  \rightarrow U^j$, for $j=1,2, \ldots, n$, are $\bigl\{\mathcal{F}_s^t\bigr\}_{s \ge t}$-adapted processes on $\bigl(\Omega, \mathcal{F}, \mathbb{P}\bigr)$ with 
\begin{align*}
\mathbb{E} \int_{s}^{T} \vert u_{\tau}^j \vert^2 d \tau < \infty,  \quad s \in [t, T];
\end{align*}
\item For any $x \in \mathbb{R}^d$, the FBSDEs in \eqref{Eq9} and \eqref{Eq15} admit a unique solution set \, $\bigl\{X_{\cdot}^{s,x;u_{\cdot}^{\neg j}}, Y_{\cdot}^{j,s,x;u_{\cdot}^{\neg j}}, Z_{\cdot}^{j,s,x;u_{\cdot}^{\neg j}}\bigr\}_{j=1}^n$ on $\bigl(\Omega, \mathcal{F}, \mathcal{F}^t, \mathbb{P}\bigr)$ and
\begin{align*}
\mathbf{Y}_{\cdot}^{s,x;u}(\omega) = \bigl(\begin{array}{cccc} Y_{\cdot}^{1,s,x;u_{\cdot}^{\neg 1}}(\omega), & Y_{\cdot}^{2,s,x;u_{\cdot}^{\neg 2}}(\omega), & \cdots, & Y_{\cdot}^{n,s,x;u_{\cdot}^{\neg n}}(\omega) \end{array}\bigr) \in K, \\
\quad \mathbb{P}-\text{\it almost}\,\, \omega \in \Omega, \quad \forall s \in [t, T]. 
\end{align*}
\end{itemize}
\end{definition}

Then, with restriction to the above admissible system, we can state the risk-averse decision problem as follows.

{\bf Problem:} Find an $n$-tuple of optimal preferable decisions for the risk-averse agents, i.e., $\hat{u}_{\cdot}=(\begin{array}{cccc} \hat{u}_{\cdot}^{1}, & \hat{u}_{\cdot}^{2}, & \cdots, & \hat{u}_{\cdot}^{n} \end{array}) \in \bigotimes\nolimits_{i=1}^n \mathcal{U}_{[t,T]}^i$, with $\xi_j^{Target} \in L^2(\Omega, \mathcal{F}_T, \mathbb{P}; \mathbb{R})$, for $j \in \{1,2, \ldots, n\}$, such that
\begin{align}
 \hat{u}_{\cdot}^{j} \in \Bigl\{ \arginf J_j &\bigr[\bigl(u^{\neg j}\bigr)\bigl] \Bigl \vert \hat{u}_{\cdot} \, \text{satisfies equation}~\eqref{Eq30} \, \text{and}\notag \\
  &  u_{\cdot}^{\neg j} = (\begin{array}{ccccccc} \hat{u}_{\cdot}^1, &\ldots, & \hat{u}_{\cdot}^{j-1}, & u_{\cdot}^j, & \hat{u}_{\cdot}^{j+1}, & \cdots, & \hat{u}_{\cdot}^n \end{array}) \in \bigotimes\nolimits_{i=1}^n \mathcal{U}_{[t,T]}^i, \quad \quad \notag \\
& \quad\quad\quad \text{with restriction to} \,\,\Sigma_{[0, T]} \Bigr\}. \label{Eq31}
\end{align}
Furthermore, the accumulated risk-costs $J_j$, for $j=1,2, \dots, n$, over the time-interval $[0, T]$ are given  
\begin{align}
J_j\bigr[\bigl(u^{\neg j}\bigr)\bigl] &= \int_0^T c_j\bigl(s, X_s^{0,x;u_{\cdot}^{\neg j}}, u_s^j\bigr) ds + \Psi_j(X_T^{0,x;u_{\cdot}^{\neg j}}), \notag \\
&\quad\quad X_0^{0,x;u_{\cdot}^{\neg j}} = x, \quad \text{and} \quad \Psi_j(X_T^{0,x;u_{\cdot}^{\neg j}}) = \xi_j^{Target}.  \label{Eq32}
\end{align}

In the following section, we establish the existence of optimal risk-averse solutions, in the sense of viscosity, for the risk-averse decision problem in \eqref{Eq31} with restriction to $\Sigma_{[0, T]}$. 

\section{Main results} \label{S4}
In this section, we present our main results, where we introduce a framework that requires a ``rational" cooperation among the risk-averse agents so as to achieve an overall optimal risk-averseness (in the sense of Pareto optimality). Moreover, such a framework allows us to establish the existence of optimal risk-averse solutions, in the sense of viscosity, to the associated risk-averse dynamic programming equations. 

\begin{proposition} \label{P1}
Suppose that the generator functional $g$ satisfies Assumption~\ref{AS1}. Further, let the statements in \eqref{Eq8} along with \eqref{Eq10} hold true. Then, for any $(t,x) \in [0, T] \times \mathbb{R}^d$ and for every $u_{\cdot}^{\neg j} = (\begin{array}{ccccccc} \hat{u}_{\cdot}^1, &\ldots, & \hat{u}_{\cdot}^{j-1}, & u_{\cdot}^j, & \hat{u}_{\cdot}^{j+1}, & \cdots, & \hat{u}_{\cdot}^n \end{array}) \in \bigotimes\nolimits_{i=1}^n \mathcal{U}_{[t,T]}^i$ and $j \in \{1,2, \ldots, n\}$, restricted to $\Sigma_{[t, T]}$, the FBSDEs in \eqref{Eq9} and \eqref{Eq15} admit unique adapted solutions
\begin{eqnarray}
\left.\begin{array}{c}
X_{\cdot}^{t,x;u_{\cdot}^{\neg j}} \in  \mathcal{S}^2\bigl(t, T; \mathbb{R}^{d} \bigr) ~~~~~~\\
\bigl(Y_{\cdot}^{j, t,x;u_{\cdot}^{\neg j}}, Z_{\cdot}^{j,t,x;u_{\cdot}^{\neg j}}\bigr) \in \mathcal{S}^2\bigl(t, T; \mathbb{R} \bigr) \times  \mathcal{H}^2\bigl(t, T; \mathbb{R}^{d} \bigr), \,\, j =1, 2, \ldots, n
\end{array}\right\}  \label{Eq33}
\end{eqnarray}
Moreover, the risk-values $V_j^{u^j}\bigl(t, x\bigr)$, for $j=1,2, \ldots, n$, are deterministic.
\end{proposition}

\begin{lemma} \label{L2}
Let $(t,x) \in [0, T] \times \mathbb{R}^d$ and $u_{\cdot}^{\neg j} = (\begin{array}{ccccccc} \hat{u}_{\cdot}^1, &\ldots, & \hat{u}_{\cdot}^{j-1}, & u_{\cdot}^j, & \hat{u}_{\cdot}^{j+1}, & \cdots, & \hat{u}_{\cdot}^n \end{array}) \\ \in \bigotimes\nolimits_{i=1}^n \mathcal{U}_{[t,T]}^i$, for $j \in \{1,2, \ldots, n \}$, be restricted to $\Sigma_{[t, T]}$. Then, for any $r \in [t, T]$ and $\mathbb{R}^d$-valued $\mathcal{F}_r^t$-measurable random variable $\eta$, we have
\begin{align}
 V_j^{u^j}\bigl(r, \eta\bigr) &= Y_r^{j, t,x;u_{\cdot}^{\neg j}} \notag\\
                                    &\triangleq \rho_{r, T}^{g_j} \Bigl[\int_r^T c_j\bigl(s, X_s^{r,\eta;u_{\cdot}^{\neg j}}, u_s^j\bigr) ds + \Psi_j(X_T^{r,\eta;u_{\cdot}^{\neg j}}) \Bigr], \quad \mathbb{P}{\text-a.s.} \label{Eq34}
\end{align}
\end{lemma}

\begin{proposition} \label{P2}
Let $u_{\cdot}^{\neg j} = (\begin{array}{ccccccc} \hat{u}_{\cdot}^1, &\ldots, & \hat{u}_{\cdot}^{j-1}, & u_{\cdot}^j, & \hat{u}_{\cdot}^{j+1}, & \cdots, & \hat{u}_{\cdot}^n \end{array}) \in \bigotimes\nolimits_{i=1}^n \mathcal{U}_{[0,T]}^i$, for $j \in \{1,2, \ldots, n \}$, be restricted to $\Sigma_{[t, T]}$. Suppose that the system of semilinear parabolic PDEs in \eqref{Eq20} enjoys the viability property w.r.t. the closed convex set $K$. Then, there exists a constant $C > 0$ such that $d_{K}^2(\cdot)$ is twice differentiable at $y$ and
\begin{align}
\bigl\langle y - \Pi_{K}(y),\, \mathbf{G}(t, &x, y, z\sigma(t,x, u^{\neg j}) \bigr\rangle \notag \\
& \le \frac{1}{4}\bigl\langle D^2(d_{K}^2(y)) z\sigma(t,x, u^{\neg j}),\,z\sigma(t,x, u^{\neg j}) \bigr\rangle + C d_{K}^2(y), \notag \\
& \quad\quad \quad \forall (t, x, y, z) \in [0, T] \times \mathbb{R}^d \times \mathbb{R}^n \times \mathcal{L}(\mathbb{R}^d; \mathbb{R}^n), \notag \\
& \quad\quad \quad\quad \forall j \in \{1,2,\ldots, n\}. \label{Eq35} 
\end{align}
\end{proposition}

Suppose that Proposition~\ref{P2} holds true, i.e., the system of semilinear parabolic PDEs in \eqref{Eq20} enjoys viability property w.r.t. the closed convex set $K$.\footnote{Here, we also assume that $K$ enjoys the BSVP for the equation in \eqref{Eq15}.} Further, assume that the set on the right-hand side of equation~\eqref{Eq30} is nonempty. Moreover, for $t \in [0,T]$ and $\hat{u}_{\cdot}=(\begin{array}{cccc} \hat{u}_{\cdot}^{1}, & \hat{u}_{\cdot}^{2}, & \cdots, & \hat{u}_{\cdot}^{n} \end{array}) \in \bigotimes\nolimits_{i=1}^n \mathcal{U}_{[t,T]}^i$, with restriction to $\Sigma_{[t, T]}$, let $\tilde{u}_{\cdot}^{\neg j}$ and $u_{\cdot}^{\neg j}$, $j \in \{1, 2, \ldots, n\}$, be two $n$-tuple decisions from $\bigotimes\nolimits_{i=1}^n \mathcal{U}_{[t,T]}^i$, i.e.,
\begin{eqnarray}
\left.\begin{array}{r}
  \tilde{u}_{\cdot}^{\neg j} = (\begin{array}{ccccccc} \hat{u}_{\cdot}^1, &\ldots, & \hat{u}_{\cdot}^{j-1}, & \tilde{u}_{\cdot}^j, & \hat{u}_{\cdot}^{j+1}, & \cdots, & \hat{u}_{\cdot}^n \end{array}) \in \bigotimes\nolimits_{i=1}^n \mathcal{U}_{[t,T]}^i \\
  u_{\cdot}^{\neg j} = (\begin{array}{ccccccc} \hat{u}_{\cdot}^1, &\ldots, & \hat{u}_{\cdot}^{j-1}, & u_{\cdot}^j, & \hat{u}_{\cdot}^{j+1}, & \cdots, & \hat{u}_{\cdot}^n \end{array}) \in \bigotimes\nolimits_{i=1}^n \mathcal{U}_{[t,T]}^i\\
   j \in \{1, 2, \ldots, n\}
\end{array}\right\}. \notag
\end{eqnarray}
Then, we can define the following partial ordering on $K$ by
\begin{align}
&\bigl(\begin{array}{cccc} \rho_{t, T}^{g_1} \bigl[\xi_{t,T}^1\bigl(\tilde{u}^{\neg 1}\bigr)\bigl], & \rho_{t, T}^{g_2} \bigl[\xi_{t,T}^2\bigl(\tilde{u}^{\neg 2}\bigr)\bigl], & \cdots, & \rho_{t, T}^{g_n} \bigl[\xi_{t,T}^n\bigl(\tilde{u}^{\neg n}\bigr)\bigl] \end{array}\bigr) \notag \\
& \quad \quad \quad \prec \bigl(\begin{array}{cccc} \rho_{t, T}^{g_1} \bigl[\xi_{t,T}^1\bigl(u^{\neg 1}\bigr)\bigl], & \rho_{t, T}^{g_2} \bigl[\xi_{t,T}^2\bigl(u^{\neg 2}\bigr)\bigl], & \cdots, & \rho_{t, T}^{g_n} \bigl[\xi_{t,T}^n\bigl(u^{\neg n}\bigr)\bigl] \end{array}\bigr), \label{Eq36}
\end{align}
if $\rho_{t, T}^{g_j} \bigl[\xi_{t,T}^j\bigl(\tilde{u}^{\neg j}\bigr)\bigl] \le \rho_{t, T}^{g_j} \bigl[\xi_{t,T}^j\bigl(u^{\neg j}\bigr)\bigl]$ for all $j=1,2, \ldots, n$, with strict inequality for at least one $j \in \{1,2, \ldots, n\}$. Furthermore, we say that
\begin{align}
\bigl(\begin{array}{cccc} \rho_{t, T}^{g_1} \bigl[\xi_{t,T}^1\bigl(\hat{u}\bigr)\bigl], & \rho_{t, T}^{g_2} \bigl[\xi_{t,T}^2\bigl(\hat{u}\bigr)\bigl], & \cdots, & \rho_{t, T}^{g_n} \bigl[\xi_{t,T}^n\bigl(\hat{u}\bigr)\bigl] \end{array}\bigr) \in K \label{Eq37}
\end{align}
is a Pareto equilibrium, in the sense of viscosity solutions, if there is no 
\begin{align}
\bigl(\begin{array}{cccc} \rho_{t, T}^{g_1} \bigl[\xi_{t,T}^1\bigl(u^{\neg 1}\bigr)\bigl], & \rho_{t, T}^{g_2} \bigl[\xi_{t,T}^2\bigl(u^{\neg 2}\bigr)\bigl], & \cdots, & \rho_{t, T}^{g_n} \bigl[\xi_{t,T}^n\bigl(u^{\neg n}\bigr)\bigl] \end{array}\bigr) \in K \label{Eq38}
\end{align}
for which
\begin{align}
&\bigl(\begin{array}{cccc} \rho_{t, T}^{g_1} \bigl[\xi_{t,T}^1\bigl(\hat{u}\bigr)\bigl], & \rho_{t, T}^{g_2} \bigl[\xi_{t,T}^2\bigl(\hat{u}\bigr)\bigl], & \cdots, & \rho_{t, T}^{g_n} \bigl[\xi_{t,T}^n\bigl(\hat{u}\bigr)\bigl] \end{array}\bigr) \notag \\
& \quad \quad \quad \prec \bigl(\begin{array}{cccc} \rho_{t, T}^{g_1} \bigl[\xi_{t,T}^1\bigl(u^{\neg 1}\bigr)\bigl], & \rho_{t, T}^{g_2} \bigl[\xi_{t,T}^2\bigl(u^{\neg 2}\bigr)\bigl], & \cdots, & \rho_{t, T}^{g_n} \bigl[\xi_{t,T}^n\bigl(u^{\neg n}\bigr)\bigl] \end{array}\bigr). \label{Eq39}
\end{align}
 
Then, with restriction to $\Sigma_{[t, T]}$, we can characterize the optimal decisions for the risk-averse agents as follows.

\begin{proposition} \label{P3}
Suppose that Proposition~\ref{P2} holds true and let $\varphi \in C_b^{1,2}([0, T] \times \mathbb{R}^d; \mathbb{R}^n)$ satisfy \eqref{Eq20} with $\varphi\bigl(T, x\bigr)=\mathbf{\Psi}(x)$ for $x \in \mathbb{R}^d$. Then, $\varphi_j\bigl(t, x\bigr) \le V_j^{u^j}\bigl(t, x\bigr)$ for $u_{\cdot}^{\neg j} = (\begin{array}{ccccccc} \hat{u}_{\cdot}^1, &\ldots, & \hat{u}_{\cdot}^{j-1}, & u_{\cdot}^j, & \hat{u}_{\cdot}^{j+1}, & \cdots, & \hat{u}_{\cdot}^n \end{array}) \in \bigotimes\nolimits_{i=1}^n \mathcal{U}_{[0,T]}^i$, for $j \in \{1,2, \ldots, n \}$, with restriction to $\Sigma_{[t, T]}$, and for all $(t,x) \in [0, T] \times \mathbb{R}^d$. Further, if an admissible optimal decision process $\hat{u}_{\cdot}^j \in \mathcal{U}_{[t,T]}^j$ exists, for almost all $(s, \Omega) \in [t, T] \times \Omega$, together with the corresponding solution $X_s^{t,x;\hat{u}_{\cdot}}$, and satisfies
\begin{align}
 \hat{u}_s^j \in &\arginf_{u_{\cdot}^j \in \mathcal{U}_{[t,T]}^j \bigl \vert \Sigma_{[t, T]}} \Bigl\{ \mathcal{L}_{s}^{u^{\neg j}} \varphi_j\bigl(s, X_s^{t,x;u_{\cdot}^{\neg j}}\bigr) \notag \\
 &\hspace{0.05in} + \,g_j\bigl(s, X_s^{t,x;u_{\cdot}^{\neg j}}, \varphi \bigl(s, X_s^{t,x;u_{\cdot}^{\neg j}}\bigr), D_x\varphi_j \bigl(s, X_s^{t,x;u_{\cdot}^{\neg j}}\bigr) \cdot \sigma\bigl(s, X_s^{t,x;w}, u_s^{\neg j}\bigr)\bigr)\Bigr\}. \label{Eq40}
\end{align}
Then, $\varphi_j\bigl(t, x\bigr) = V_j^{\hat{u}^j}\bigl(t, x\bigr)$ for $j \in \{1,2, \ldots, n \}$ and for all $(t,x) \in [0, T] \times \mathbb{R}^d$. Moreover, corresponding to the $n$-tuple of optimal risk-averse decisions $\hat{u}_{\cdot} \in \bigotimes\nolimits_{i=1}^n \mathcal{U}_{[t,T]}^i$, there exists a Pareto equilibrium
\begin{align}
\bigl(\begin{array}{cccc} \rho_{t, T}^{g_1} \bigl[\xi_{t,T}^1\bigl(\hat{u}\bigr)\bigl], & \rho_{t, T}^{g_2} \bigl[\xi_{t,T}^2\bigl(\hat{u}\bigr)\bigl], & \cdots, & \rho_{t, T}^{g_n} \bigl[\xi_{t,T}^n\bigl(\hat{u}\bigr)\bigl] \end{array}\bigr) \in K \label{Eq41}
\end{align}
such that 
\begin{align}
&\bigl(\begin{array}{cccc} \rho_{t, T}^{g_1} \bigl[\xi_{t,T}^1\bigl(\hat{u}\bigr)\bigl], & \rho_{t, T}^{g_2} \bigl[\xi_{t,T}^2\bigl(\hat{u}\bigr)\bigl], & \cdots, & \rho_{t, T}^{g_n} \bigl[\xi_{t,T}^n\bigl(\hat{u}\bigr)\bigl] \end{array}\bigr) \notag \\
& \quad \quad \quad \prec \bigl(\begin{array}{cccc} \rho_{t, T}^{g_1} \bigl[\xi_{t,T}^1\bigl(u^{\neg 1}\bigr)\bigl], & \rho_{t, T}^{g_2} \bigl[\xi_{t,T}^2\bigl(u^{\neg 2}\bigr)\bigl], & \cdots, & \rho_{t, T}^{g_n} \bigl[\xi_{t,T}^n\bigl(u^{\neg n}\bigr)\bigl] \end{array}\bigr) \,\, \text{on} \,\, K, \label{Eq42}
\end{align}
for all $t \in [0, T]$ and for $\bigl\{\xi_j^{Target}\bigr\}_{j=1}^n$ from $L^2(\Omega, \mathcal{F}_T, \mathbb{P}$.
\end{proposition}

\section{Further remarks} \label{S5}

In this section, we briefly comment on the problem formulation, where the risk-averse decision framework of Section~\ref{S3} -- in which results from the dynamic risk analysis implicitly constitute part of the information used in the context of the risk-averse criteria -- requires each of the risk-averse agents to respond optimally (in the sense of {\it best-response correspondence}) to the decisions of the other risk-averse agents. 

Note that, for example, see equations~\eqref{Eq36}--\eqref{Eq38} for the notion of Pareto equilibrium, and see also equation~\eqref{Eq40} for consistent optimal decisions that are all well defined concepts in the context of risk-aversion problem in \eqref{Eq31}, with accumulated risk-costs of \eqref{Eq7}. Here, we remark that, for every $\bigl\{\xi_j^{Target}\bigr\}_{j=1}^n$ from $L^2(\Omega, \mathcal{F}_T, \mathbb{P}; \mathbb{R})$ and for all $t \in [0, T]$, if there exists an $n$-tuple of optimal risk-averse decisions, i.e., $\hat{u}_{\cdot} \in \bigotimes\nolimits_{i=1}^n \mathcal{U}_{[t,T]}^i$, such that, for any $x \in \mathbb{R}^d$, the FBSDEs in \eqref{Eq9} and \eqref{Eq15} (cf. equations~\eqref{Eq13} and \eqref{Eq14}) admit a unique solution set $\bigl\{X_{\cdot}^{t,x;u_{\cdot}^{\neg j}}, Y_{\cdot}^{j,t,x;u_{\cdot}^{\neg j}}, Z_{\cdot}^{j,t,x;u_{\cdot}^{\neg j}}\bigr\}_{j=1}^n$ on $\bigl(\Omega, \mathcal{F}, \mathcal{F}^t, \mathbb{P}\bigr)$ and
\begin{align*}
\mathbf{Y}_{\cdot}^{s,x;u}(\omega) = \bigl(\begin{array}{cccc} Y_{\cdot}^{1,s,x;u_{\cdot}^{\neg 1}}(\omega), & Y_{\cdot}^{2,s,x;u_{\cdot}^{\neg 2}}(\omega), & \cdots, & Y_{\cdot}^{n,s,x;u_{\cdot}^{\neg n}}(\omega) \end{array}\bigr) \in K, \\
\quad \mathbb{P}-\text{\it almost}\,\, \omega \in \Omega, \quad \forall s \in [t, T]. 
\end{align*}
Then, verifying the above condition amounted to solving the stochastic target problem, which can be specified by a set of all acceptable risk-exposures, when $t=0$, vis-\'{a}-vis some uncertain future costs or outcomes specified by a set of random variables $\bigl\{\xi_j^{Target}\bigr\}_{j=1}^n$ from $L^2(\Omega, \mathcal{F}_T, \mathbb{P}; \mathbb{R})$. 

On other hand, assume that the exact information about $\xi_j^{Target} \in L^2(\Omega, \mathcal{F}_T, \mathbb{P}; \mathbb{R})$, for $j=1,2,\ldots, n$, are not a-priorly known, but we know that such information can be obtained from the following allocation
\begin{align*}
L = \sum\nolimits_{j=1}^n \alpha_j \xi_j^{Target} \in L^2(\Omega, \mathcal{F}_T, \mathbb{P}; \mathbb{R}),
\end{align*}
where $L \in L^2(\Omega, \mathcal{F}_T, \mathbb{P}; \mathbb{R})$ is assumed known and $\sum\nolimits_{j=1}^n \alpha_j =1$, for some $\alpha_j \ge 0$,  $j=1,2,\ldots, n$. Furthermore, if there exists an $n$-tuple of optimal decisions, i.e., $\hat{u}_{\cdot} \in \bigotimes\nolimits_{i=1}^n \mathcal{U}_{[t,T]}^i$, for the risk-averse agents, then we can define the set of optimally allocated risk-exposures as follows
\begin{align*}
\mathcal{A}_{0}(L) = \Biggl\{\bigl(\begin{array}{cccc} \rho_{0, T}^{g_1} \bigl[\xi_{0,T}^1\bigl(\hat{u}\bigr)\bigl], & \rho_{0, T}^{g_2} \bigl[\xi_{0,T}^2\bigl(\hat{u}\bigr)\bigl], & \cdots, & \rho_{0, T}^{g_n} \bigl[\xi_{0,T}^n\bigl(\hat{u}\bigr)\bigl] \end{array}\bigr) \in K  \, \Bigl\vert \, \quad &\\
\,\,L = \sum\nolimits_{j=1}^n \alpha_j \xi_j^{Target} \in L^2(\Omega, \mathcal{F}_T, \mathbb{P}; \mathbb{R}) \,\,\,\text{and} \,\,\, \sum\nolimits_{j=1}^n \alpha_j =1, \,\, & \\
\text{with} \,\,\, \alpha_j \ge 0,\,\, j=1,2,\ldots, n & \Biggr\},
\end{align*}
which provides further useful information to characterize all Pareto equilibria w.r.t. the risk-averse agents.

\end{document}